%Paper: q-alg/9502011
%From: yamada@math.metro-u.ac.jp
%Date: Sat, 18 Feb 95 17:08:18 JST

\input amstex
\documentstyle{amsppt}

\magnification = 1200
\hsize = 16truecm
\vsize = 21truecm
\parskip = 6pt

\def\N{\Bbb N}
\def\C{\Bbb C}
\def\Z{\Bbb Z}
\def\A{A^{(1)}_1}
\def\g{\goth{g}}
\def\Sy{\goth{S}}
\def\p{\partial}
\def\a{\alpha}

\vskip 1cm

\centerline{\bf WEIGHT VECTORS OF THE BASIC $\A$-MODULE}

\centerline{\bf AND THE LITTLEWOOD-RICHARDSON RULE}

\rightheadtext{WEIGHT VECTORS OF THE BASIC $\A$-MODULE}
\leftheadtext{S. ARIKI, T. NAKAJIMA AND H.-F. YAMADA}

\vskip 1.0in

\centerline{  Susumu ARIKI  }
\vskip 0.1in
\centerline{  Division of Mathematics,  }
\centerline{  Tokyo University of Mercantile Marine,  }
\centerline{  2-1-6, Etchujima,  }
\centerline{  Koto-ku, Tokyo 135, Japan }
\vskip 0.2in
\centerline{  Tatsuhiro NAKAJIMA  }
\vskip 0.1in
\centerline{  Department of Physics, }
\centerline{  Faculty of Science,  }
\centerline{  Tokyo Metropolitan University,  }
\centerline{  1-1, Minami-Ohsawa, Hachioji-shi,  }
\centerline{  Tokyo 192-03, Japan   }
\vskip 0.2in
\centerline{  Hiro-Fumi YAMADA  }
\vskip 0.1in
\centerline{
Department of Mathematics,  }
\centerline{  Faculty of Science,  }
\centerline{  Tokyo Metropolitan University,  }
\centerline{  1-1, Minami-Ohsawa, Hachioji-shi,  }
\centerline{  Tokyo 192-03, Japan   }
\vskip 1.0in

\centerline{\bf Abstract}

The basic representation of $\A$ is studied. The weight
vectors are represented in terms of Schur functions. A suitable
base of any weight space is given. Littlewood-Richardson rule
appears in the linear relations among weight vectors.

\pagebreak
\document
\heading
\S 0 INTRODUCTION
\endheading

The aim of this letter is to give an explicit expression of the
weight vectors of the basic $\A$-module realized on the polynomial
ring of infinitely many variables.

In 1978, Lepowsky and Wilson [5] constructed the basic representation
of the affine Lie algebra $\A$ by making use of the vertex operator.
This construction was generalized to other types of affine Lie
algebras [4] and applied to a study of nonlinear integrable
differential equations such as  KP, BKP and KdV hierarchies [1].
Roughly speaking, weighted homogeneous polynomial solutions
($\tau$-functions) of these hierarchies are weight vectors
whose weights lie on the Weyl group orbit through the highest weight,
namely the maximal weights.
They are expressed by means of the Schur functions or the
Schur $Q$-functions, reflecting that the formal solutions constitute
an infinite dimensional Grassmann manifold or its submanifold.

In this letter we show that the weight space of any weight of the
basic $\A$-module is spanned by $2$-reduced Schur functions. We will
choose a suitable base of any weight space and discuss the linear
relations among the weight vectors. To our surprise, the
Littlewood-Richardson rule [2,6] appears in the linear relations.
The formula obtained can also be viewed as an identity for $2$-modular
characters of the symmetric group. We believe that there is a deep
connection between affine Lie algebras and modular representations
of the symmetric group.

We also believe that the Virasoro algebra is related with this
story. Wakimoto [9,10] proved that the sum of weight spaces of weights
$\{\Lambda-n\delta; n\in \N \}$ admits a Virasoro action,
where $\Lambda$ is a maximal weight and $\delta$ is the
fundamental imaginary root.
It is  irreducible in the case of the basic $\A$-module.
As a corollary, the Schur functions indexed by staircase Young diagrams are
characterized as singular vectors of the Virasoro representation.
The relation between this Virasoro representation and our weight
vectors must be clarified. ({\it See also} [7,11].)

\heading
\S 1 REVIEW OF THE BASIC $\A$-MODULE
\endheading

We first review some ingredients of a realization
of the basic $\A$-module [3,5].
Let $\g=\A$ be the affine Lie algebra corresponding to the
Cartan matrix $\left(\matrix 2 & -2 \\ -2 & 2 \\ \endmatrix \right)$
with the standard Chevalley generators
$\{e_0, e_1, f_0, f_1, \alpha_0^{\vee}, \alpha_1^{\vee}\}$.
The basic representation of $\g$ is realized on the space
of polynomials of infinitely many variables
$$V=\C[t_1,t_3,t_5,\ldots]$$
as follows. For any odd natural number $j$, let
$a_j=\frac{\partial}{\partial t_j}$ and $a_{-j}=jt_j$. Then
$\{a_j (j \in \Z, odd), Id\}$ span the infinite dimensional
Heisenberg algebra acting on $V$, which is a subalgebra of $\g$.
The action of $\g$ is constructed by the so called vertex
operator. Let $p$ be an indeterminate and put
$$\align \xi(t,p)&=\sum_{j\geq1,odd} t_jp^j
                   =\sum_{j\geq1,odd}\frac{a_{-j}}{j}p^j, \\
         \xi(\tilde \p,p^{-1})&=
         \sum_{j\geq1,odd} \frac{1}{j} \frac{\p}{\p t_j} p^{-j}
                  =\sum_{j\geq1,odd} \frac{a_j}{j} p^{-j}.
 \endalign$$
The vertex operator is defined by
$$ X(p)=-\frac{1}{2} e^{2\xi(t,p)}e^{-2\xi(\tilde \p,p^{-1})}. $$
Expanding $X(p)$ as a formal power series of $p$ and $p^{-1}$:
$$ X(p)=\sum_{k\in \Z} X_k p^{-k},$$
we have differential operators $X_k$ ($k\in \Z$) acting on $V$.
It is proved [5] that operators $a_j$ ($j\in \Z$,odd), $X_k$ (
$k\in \Z$) and identity constitute the affine Lie algebra
$\g=\A$, namely the basic representation of $\g$. This is the irreducible
highest weight $\g$-module with highest weight $\Lambda_0$, where
$\Lambda_0(\a_0^{\vee})=1$, $\Lambda_0(\a_1^{\vee})=0$.
Let $\a_0$ and $\a_1$ be the simple roots of $\g$,
and $\delta=\a_0+\a_1$ be the fundamental imaginary root.
It is well-known [3] that the set of weights $P$ of the basic
$\g$-module is given by
$$
P=\left\{\Lambda_0+q \delta+p \alpha_1; p,q\in \Z, q\leq-p^2 \right\}.
$$
A weight $\Lambda$ on the parabola $q=-p^2$ is
said to be maximal in the sense that $\Lambda+\delta$
is no longer a weight. Maximal weights consist a single
Weyl group orbit.
For a maximal weight $\Lambda$ the weight vector is expressed by
the Schur function. For any Young diagram $Y$ of $N$ cells,
the Schur function indexed by $Y$ is defined by
$$ S_Y(t)=\sum_{\nu_1+2\nu_2+\cdots=N} \chi_Y(\nu)
          \frac{t_1^{\nu_1}t_2^{\nu_2}\cdots}{\nu_1!\nu_2!\cdots}, $$
where $\chi_Y(\nu)$ is the character value of the irreducible
representation $Y$ of the symmetric group $\Sy_N$, evaluated at the
conjugacy class of the cycle type
$\nu=(1^{\nu_1}2^{\nu_2}\cdots N^{\nu_N})$ [6].
The Schur function $S_Y(t)$ is obviously a weighted homogeneous
(deg$t_j=j$) polynomial of degree $|Y|$.
For a non-negative integer $r$ put $K_r=(r,r-1,r-2,\ldots,2,1)$ be
the staircase Young diagram of length $r$. In terminology of
modular representations of the symmetric group, $K_r$ are called
$2$-cores since they do not have $2$-hooks [8].
We remark that the Schur functions $S_{K_r}(t)$ ($r=0,1,2,\ldots$)
do not depend on $t_{2j}$ ($j=1,2,\ldots$), namely elements of  $V$,
because of the Murnaghan-Nakayama formula [2]. The maximal weight vectors
are $S_{K_r}(t)$ for $r=0,1,2,\ldots$ [1]. We denote by $\Lambda_r$
the maximal weight whose weight vector is $S_{K_r}(t)$. According
to the theory of hierarchies of nonlinear integrable systems, these
maximal weight vectors exhaust the weighted homogeneous polynomial
$\tau$-functions  of the KdV hierarchy [1].

The subspace of $V$ consisting of weighted homogeneous polynomials
of degree $n$ has dimension $p^{odd}(n)$, the number of partitions
of $n$ into odd positive integers. If we put $\phi(q)=\prod_{j\geq 1}
(1-q^j)$, the generating function is
$$\sum_{n=0}^{\infty} p^{odd}(n)q^n=\frac{\phi(q^2)}{\phi(q)}.$$

Let $\mu(n)$ be the weight multiplicity of the weight $\Lambda-n\delta$
for a maximal weight $\Lambda$, which is independent of the choice
of the maximal weight, since the Weyl group preserves the weight
multiplicity. It is obvious that degree of a maximal weight
vector is of the form $2m^2+m$ ($m\in \Z$). Therefore, by using an
identity of Gauss [3,p241], we see that
$$\frac{\phi(q^2)}{\phi(q)}=\sum_{m\in \Z}\left(
  \sum_{n=0}^{\infty} \mu(n)q^{2n} \right) q^{2m^2+m}
  =\frac{\phi(q^2)^2}{\phi(q)} \sum_{n=0}^{\infty}\mu(n)q^{2n}, $$
and hence $\mu(n)=p(n)$, the number of partition of $n$ into positive
integers.

\heading
\S 2 BASES FOR WEIGHT SPACES
\endheading

We define the $2$-quotient for a given Young diagram [2,8]. Let
$Y=(y_1,\ldots,y_n)$ \linebreak $(y_1 \geq \cdots \geq y_n \geq 0)$ be a
Young diagram. We always assume $n$ to be even. Consider the ``
Maya diagram'' or the ``$\beta$-set'' $X=(x_1,\ldots,x_n)$ where
$x_j=y_j+(n-j)$ for $1\leq j \leq n$. For $i=0, 1$ let
$$
X^{(i)}=\left\{\xi^{(i)}\in \N; 2\xi^{(i)}+i=x_j \text{ for some } j \right\}.
$$
If we have $X^{(i)}=\left\{\xi^{(i)}_1,\ldots,\xi^{(i)}_{m^{(i)}}\right\} \,
(\xi^{(i)}_1>\cdots>\xi^{(i)}_{m^{(i)}}\geq 0)$, then we define
the Young diagram $Y^{(i)}$ by
$$
Y^{(i)}=\left(\xi^{(i)}_1-(m^{(i)}-1), \xi^{(i)}_2-(m^{(i)}-2),\ldots
\xi^{(i)}_{m^{(i)}}\right).
$$
The pair $(Y^{(0)},Y^{(1)})$ of Young diagrams is called the $2$-quotient
of $Y$. The $2$-core of $Y$ is described as follows. If
$|X^{(0)}|-|X^{(1)}|=r \geq 1$, then the $2$-core of $Y$ is
$K_{r-1}$ and if $|X^{(1)}|-|X^{(0)}|=r \geq 0$, then it is $K_r$.
In this fashion we can attach a triplet $(K; Y^{(0)}, Y^{(1)})$
of Young diagrams
for any Young diagram $Y$, where $K$ is the $2$-core of $Y$ and
$(Y^{(0)},Y^{(1)})$ is the $2$-quotient of $Y$. It is easily shown
that this correspondence is one-to-one and $|Y|=2(|Y^{(0)}|+|Y^{(1)}|)
+|K|$. Denote by $\tau(Y)$ the triplet $(K; Y^{(0)},Y^{(1)})$
corresponding to $Y$.

\proclaim{Example}
$$
\alignat 2
&Y=(4,3,1^2), & \qquad &X=(7,5,2,1), \\
&X^{(0)}=(1), & \qquad &X^{(1)}=(3,2,0), \\
&Y^{(0)}=(1), & \qquad &Y^{(1)}=(1^2,0), \\
&K=K_2=(2,1).
\endalignat
$$
\endproclaim

We now describe the weight vectors of the basic $\A$-module by
means of the Schur functions. To this end we denote by
$S^{red}_Y(t) \in V$ the $2$-reduced Schur function indexed
by $Y$, which is by definition,
$$
S^{red}_Y(t)=S_Y(t)|_{t_2=t_4=\cdots=0}.
$$
By using the Boson-Fermion correspondence established by Date et al.[1],
we can see the following.

\proclaim{Proposition 1} The $2$-reduced Schur function $S^{red}_Y(t)$
is a weight vector of weight $\Lambda_r-n\delta$ if
$\tau(Y)=(K_r; Y^{(0)},Y^{(1)})$ with $|Y^{(0)}|+|Y^{(1)}|=n$.
\endproclaim

As we have seen in \S 1, the multiplicity of each weight is
expressed by the number of partitions for any maximal weight
$\Lambda_r$, {\it i.e.}, $mult(\Lambda_r-n\delta)=p(n)$.
Hence the weight vectors found above satisfy linear relations in general.
For example, the Young diagrams $(4), (3,1), (2^2), (2,1^2), (1^4)$
determine the $2$-reduced Schur functions of the same weight
$\Lambda_0-2\delta$ whose multiplicity equals $p(2)=2$. It is easily
checked that
$$
\left\{
\aligned
S^{red}_{(4)}(t)&=S^{red}_{(1^4)}(t)=\frac{1}{24}t_1^4+t_1t_3, \\
S^{red}_{(3,1)}(t)&=S^{red}_{(2,1^2)}(t)=\frac{1}{8}t_1^4, \\
S^{red}_{(2^2)}(t)&=S^{red}_{(2,1^2)}(t)-S^{red}_{(4)}(t)
                   =\frac{1}{12}t_1^4-t_1t_3.
\endaligned
\right.
$$
Therefore the next problem is to find a suitable base for each
weight space. The following theorem gives an answer.

\proclaim{Theorem 2} The $2$-reduced Schur functions
$$
\left\{S^{red}_Y(t); \tau(Y)=(K_r; \phi, Y^{(1)}) \text{ with }
|Y^{(1)}|=n \right\}
$$
are linearly independent and hence constitute a base for the weight space of
weight $\Lambda_r-n\delta$.
\endproclaim

Any weight vector $S^{red}_Y(t)$ can be expressed uniquely as a
linear combination of the base vectors obtained above. We now focus
on the coefficients of these expressions. Suppose that the Young
diagram $Y$ corresponds to the triplet $\tau(Y)=(K; Y^{(0)}, Y^{(1)})$.
The $2$-sign $\delta_2(K; Y^{(0)},Y^{(1)})$ is defined as follows.
If the $2$-core $K$ is obtained from $Y$ by removing a sequence of
$2$-hooks, where $q$ of them are column $2$-hooks and the others
are row $2$-hooks, then
$$
\delta_2(K; Y^{(0)}, Y^{(1)})=(-1)^q.
$$
It can be proved that $\delta_2(K; Y^{(0)}, Y^{(1)})$ does not depend
on the choice of $2$-hooks being removed. The following is our main
result.

\proclaim{Theorem 3} For such a Young diagram $Y$ that
$\tau(Y)=(K; Y^{(0)}, Y^{(1)})$ with $n=|Y^{(0)}|+|Y^{(1)}|$,
we have
$$
S^{red}_Y(t)=(-1)^{|Y^{(0)}|} \delta_2(K; Y^{(0)}, Y^{(1)})
\sum_{Z^{(1)}} LR^{Z^{(1)}}_{Y^{(0)\prime} Y^{(1)}}
\delta_2(K; \phi, Z^{(1)}) S^{red}_Z(t),
$$
where the summation runs over all Young diagrams $Z^{(1)}$ of size $n$,
the Young diagram $Z$ corresponds to $(K; \phi, Z^{(1)})$ and
$Y^{(0)\prime}$ denotes the transpose of $Y^{(0)}$.
We also denote by $LR$ the Littlewood-Richardson coefficient.
\endproclaim
The proof of this theorem is performed by analysing the Schur functions,
apart from the affine Lie algebra $\A$. Details will be published
elsewhere.

\enddocument